\documentclass[11pt,a4paper]{article}

\usepackage{amsfonts}
\usepackage{color}
\usepackage{graphicx}
\usepackage{float}
\usepackage{alltt}

\usepackage{amsthm}
\usepackage{amsmath}
\usepackage{thmtools}
\usepackage{graphics}
\usepackage{epsfig}
\usepackage{amssymb}
\usepackage{color}

\usepackage{graphicx}
\usepackage{newtxtext,newtxmath}
\usepackage{enumerate}
\usepackage{footnote}
\usepackage[colorlinks=true, allcolors=blue]{hyperref}
\usepackage[normalem]{ulem}

\usepackage[normalem]{ulem}

\usepackage{geometry}

\newcommand{\N}{\mathbb{N}}
\newcommand{\R}{\mathbb{R}}
\newcommand{\C}{\mathbb{C}}
\newcommand{\Z}{\mathbb{Z}}
\DeclareMathOperator{\diag}{diag}

\newtheorem{theorem}{Theorem}
\theoremstyle{definition} 
\newtheorem{example}{Example}
\newtheorem{proposition}{Proposition}
\newtheorem{remark}{Remark}
\newtheorem{definition}{Definition}

\title{\textbf{On the stability of $\theta$-methods for DDEs and PDDEs} }

\author{Alejandro Rodríguez-Fernández$^{(a)}$
        Jesús Martín-Vaquero$^{(b)}$
\\
           $^{(a)}$   Department of Applied Mathematics, University of Salamanca, Salamanca, Spain\\
              email: alexRF@usal.es, arf262001@gmail.com
\\
            $^{(b)}$  Department of Applied Mathematics, IUFFyM,University of Salamanca, Salamanca, Spain\\
              email: jesmarva@usal.es
}

\date{}

\begin{document}

\maketitle

\begin{abstract}
In this paper, the stability of $\theta$-methods for delay differential equations is studied based on the test equation $y'(t)=-A y(t) + B y(t-\tau)$, where $\tau$ is a constant delay and $A$ is a positive definite matrix. It is mainly considered the case where the matrices $A$ and $B$ are not simultaneosly diagonalizable and the concept of field of values is used to prove a sufficient condition for unconditional stability of these methods and another condition which also guarantees their stability, but according to the step size.
The results obtained are also simplified for the case where the matrices $A$ and $B$ are simultaneously diagonalizable and compared with other similar works for the general case. Several numerical examples in which the theory discussed here is applied to parabolic problems given by partial delay differential equations with a diffusion term and a delayed term are presented, too.


\end{abstract}

\section{Introduction}
\label{sec:1}

 Many real-life phenomena in bioscience require differential equations that depend partially on the past, not only the current state.
Examples appear in population dynamics, infectious diseases (accounting for the incubation periods), chemical and enzyme kinetics and more general control problems (see \cite{BoRi00,macdonald2008biological,MaGl77,TaMaBe00} and references therein).

They include several types of functional differential equations: delay differential equations (DDEs), neutral delay differential equations (NDDEs) \cite{Jackiewicz1987/88},
integro-differential \cite{Baker06}, partial delay differential equations (PDDEs)  \cite{KotoRK,KotoMultistep},
or stochastic delay differential equations (SDDEs) \cite{baker_buckwar_2000,KUCHLER2000189}.

However, the number of occasions on which an exact solution to such problems can be found is very limited.
Because of this, numerical methods are usually used to obtain an approximate solution.  In the scientific literature, to study the stability of different numerical methods for delay differential equations it has usually been followed the classical theory based on studying the behavior of the method when applied to the differential equation with constant delay
\begin{equation}
   \label{eq:1:EDR_EstabilidadTipica}
   y'(t) = \lambda y(t) + \gamma y(t-\tau),
\end{equation}
where $\lambda, \gamma \in \C$ (see \cite{KotoRK,KotoMultistep,Calvo,Rihan}).

This analysis is similar to the most frequent one to study other types of numerical methods (including IMEX, ETD, exponential fitting or additive semi-implicit) for differential equations and partial differential equations without delay \cite{DAUTILIA20202067,Cox2002ExponentialTD,VIGOAGUIAR200780,VAMVWa07,Higueras_2015} and references therein.
However, in this way the study of stability of the numerical methods applied to DDEs is reduced to the delay differential equations of the form
\begin{equation*}
    y'(t) = A \cdot y(t) + B \cdot y(t-\tau),
\end{equation*}
where $\tau$ is as before and $A$ and $B$ are two simultaneously diagonalizable matrices in $\C$. It would consist of studying the stability for several DDEs as (\ref{eq:1:EDR_EstabilidadTipica}), where, for each of the DDEs, the coefficient $\lambda$ would be an eigenvalue of $A$, while $\gamma$ would be the corresponding eigenvalue of $B$ with the same eigenvector as $\lambda$.

For this reason, in this paper (Sect. \ref{sec:3}) we will follow a newer theory that has only been used to study the stability of some numerical methods for problems without delay (see \cite{SeiboldTheory,SeiboldPractice}), proving that it can also be used when working with delays, and, as a result, we will be able to know the asymptotic behavior of the $\theta$-methods defined below when applied to differential equations with a constant delay such as
\begin{equation}
    y'(t)=-A y(t) + B y(t-\tau),
    \label{eq:1:EDR_EstabilidadThetaMétodo}
\end{equation}
where $A$ is now a positive definite matrix, but need not be simultaneously diagonalizable with $B$. Later, we will apply the results obtained to simplify the study of stability for equation (\ref{eq:1:EDR_EstabilidadThetaMétodo}) when the matrices A and B are simultaneously diagonalizable and, in this way, we shall compare the new results with others obtained with different procedures \cite{Calvo,Rihan,kj1994stability}, demonstrating that the new strategy can be useful to explain stability in a large number of problems and numerical methods.
Finally, although our study is based on the DDE (\ref{eq:1:EDR_EstabilidadThetaMétodo}), in Section \ref{sec:5} we will show that the theory presented here can also be used to know the stability of the aforementioned methods when we use them to calculate a numerical solution of certain partial delay differential equations that depend on both time and space.
For it, we will transform the corresponding parabolic problem into a system of only time-dependent DDEs using the method of lines (MOL), which consists of discretizing space by approximating each of the partial derivatives with respect to space by a difference equation. We will explain this approach in more detail by means of several numerical examples.

\section{Preliminaries}
\label{sec:2}

In this section, we recall some notations, definitions and preliminary materials.

Let $M \in M_N(\C)$ be a square matrix. The field of values of $M$, denoted as $F(M)$, is defined by
\begin{equation*}
    F(M) := \{ x^* Mx : x \in \C^N, \: x^*x = 1 \},
\end{equation*}
where $x^*$ denotes the conjugate transpose of a vector $x$. Moreover, we denote the spectrum of $M$ by $\sigma(M)$.

Then, for two arbitraries matrices $A,B \in M_N (\C)$, the following properties are verified:
\begin{enumerate}[(a)]
    \item $F(A)$ is convex, closed and bounded.
    \item $\sigma(A) \subseteq F(A)$.
    \item $F(A+B) \subseteq F(A) + F(B)$.
    \item $F(\alpha A) = \alpha F(A)$, $\alpha \in \C$.
    \item If the matrix $B$ is positive definite, $\sigma(A B) \subseteq F(A) F(B)$.
    \item If $A$ is normal, $F(A) = Co(\sigma(A))$, where $Co(\cdot)$ denote the closed convex hull of a set.
    \item If $A$ is hermitian, $F(A)$ is a segment of the real line.
\end{enumerate}

The reader interested in the proofs of the above properties, as well as a deeper review of the field of values of a matrix, is referred to \cite{HornJohnson,Johnson}, where it can also be found an algorithm to delimit this set.

On the other hand, given the ordinary differential equation
\begin{equation}
    y'(t) = f(t,y(t)) + g(t,y(t)),
    \label{eq:2:ODE_ThetaMétodo}
\end{equation}
we define a $\theta$-method for (\ref{eq:2:ODE_ThetaMétodo}) by
\begin{equation*}
    y_{n+1} = y_n + h (1-\theta)[f(t_n, y_n) + g(t_n, y_n)] + h \theta [f(t_{n+1}, y_{n+1}) + g(t_{n+1}, y_{n+1})],
\end{equation*}
where the parameter $\theta \in [0,1]$.

Therefore, setting $\tau=(m-u)h$, with $m \in \N$ and $u \in [0,1)$, and using linear interpolation to approximate the terms with delay, we can define a $\theta$-method for the delay differential equation
\begin{equation*}
    y'(t) = f(t,y(t)) + g(t,y(t),y(t-\tau)),
\end{equation*}
 by means of
\begin{equation}
    \begin{split}
        y_{n+1} = y_n &+ h(1-\theta)[f(t_n, y_n) + g(t_n, y_n, (1-u) y_{n-m} + u \cdot y_{n-m+1})] \\ &+h\theta[f(t_{n+1}, y_{n+1}) + g(t_{n+1}, y_{n+1}, (1-u) y_{n-m+1} + u \cdot y_{n-m+2})],
    \end{split}
    \label{eq:2:thetaMétodoEDRs}
\end{equation}
where $\theta$ is as before.

As to the order of these methods, it is verified that the error of the linear interpolation considered does not influence them (see \cite{Calvo}), so we can easily calculate it considering their linear interpolation error $T_{n+1}$ for $u=0$.
This way, we have that
\begin{equation*}
        T_{n+1}= \left (\frac{1}{2} -\theta \right) \bigg[\digamma \cdot (f(t_n,y(t_n)) +g(t_n,y(t_n),y(t_n))) +\frac{\partial f}{\partial x_1}(t_n,y(t_n))
        +\frac{\partial g}{\partial x_1}(t_n,y(t_n),y(t_n)) \bigg] h+O\left(h^2\right),
\end{equation*}
where $\frac{\partial f}{\partial x_i}$ and $\frac{\partial g}{\partial x_i}$, $i=1,2,3$, denote the partial derivative with respect to the i-th variable of $f$ and $g$, respectively, and
\begin{equation*}
    \digamma := \frac{\partial f}{\partial x_2}(t_n,y(t_n))+\frac{\partial g}{\partial x_3}(t_n,y(t_n),y(t_n))+\frac{\partial g}{\partial x_2}(t_n,y(t_n),y(t_n)).
\end{equation*}
That is, if $\theta \neq 1/2$ the method (\ref{eq:2:thetaMétodoEDRs}) has order 1. Similarly, it can also be seen that $\theta$-method with $\theta = 1/2$ has order 2.

\section{Stability}
\label{sec:3}

As indicated at the beginning of this paper, to study the stability properties of the proposed $\theta$-methods, we first consider the test problem (\ref{eq:1:EDR_EstabilidadThetaMétodo}).

A $\theta$-method for this problem is defined by
\begin{equation}
    \begin{split}
       y_{n+1} = y_n &+ h(1-\theta)[-A \cdot y_n + B \cdot ((1-u) y_{n-m} + u \cdot y_{n-m+1})] \\ &+h\theta[ -A \cdot y_{n+1} + B \cdot ((1-u) y_{n-m+1} + u \cdot y_{n-m+2})].
    \end{split}
    \label{eq:3:thetaMetodoEstabilidad1}
\end{equation}
Then, by means of $\boldsymbol{Y_n}=(y_n, y_{n-1},\ldots,y_{n-m})$, we can rewrite (\ref{eq:3:thetaMetodoEstabilidad1}) in the form
\begin{equation}
     \boldsymbol{Y_{n+1}} = W \cdot \boldsymbol{Y_n},
     \label{eq:3:thetaMetodoEstabilidad2}
\end{equation}
where
\begin{equation*}
    W =
    \begin{pmatrix}
        I-\theta h (-A) & 0 & 0 & \ldots & 0\\
        0 & I & 0 & \ldots & 0\\
        0 & 0 & I & \ldots & 0\\
        \vdots & \vdots & \vdots & \ddots & \vdots\\
        0 & 0 & 0 & \ldots & I\\
    \end{pmatrix}^{-1}
    \cdot
    \begin{pmatrix}
        I+(1-\theta) h (-A) & 0  & \ldots & 0 & h B \theta u & \Upsilon_1 &  \Upsilon_2\\
        I & 0 & \ldots & 0 & 0 & 0 & 0 \\
        0 & I & \ldots & 0 & 0 & 0 & 0 \\
        \vdots & \vdots & \ddots & \vdots & \vdots & \vdots & \vdots\\
        0 & 0 & \ldots & 0 & I & 0 & 0\\
        0 & 0 & \ldots & 0 & 0 & I & 0\\
    \end{pmatrix}
\end{equation*}
for $m \geq 3$, $\Upsilon_1 := h B (\theta (1-u) + (1-\theta) u)$, $\Upsilon_2 :=h B (1-\theta)(1-u)$ and $I$ denote the $N \times N$ identity matrix (being $N$ the dimension of the matrices $A$ and $B$).

Let $V \neq 0$ a eigenvector of $W$ with eigenvalue $\xi$. Then, due to the form of the matrix $W$, it is verified that
\begin{equation*}
    V= \big (\xi ^m v,\xi ^{m-1} v,\ldots,\xi v,v \big )
\end{equation*}
for some $v \neq 0$, $v \in \C^N$, such that
\begin{equation*}
    T(\xi) v =  0,
\end{equation*}
where
\begin{equation*}
    \begin{split}
      &T(z) :=  \frac{1}{h} a(z) I - c(z)(-A) - b(z)B,\\
      &a(z):= z^{m+1} - z^m,\\
      &b(z) := \theta (u z^2 + (1-u)z)+ (1-\theta)(u z + (1-u)),\\
      &c(z) := \theta z^{m+1} + (1-\theta) z^m.
    \end{split}
\end{equation*}

Hence, using the characteristic equation of (\ref{eq:3:thetaMetodoEstabilidad2}), we can characterize the stability of (\ref{eq:3:thetaMetodoEstabilidad1}) as follows.

\begin{proposition}
    If the matrix $T(z)$ is nonsingular for all $|z| \geq 1$, then the corresponding $\theta$-method is stable.
\end{proposition}

On the other hand, if we multiply $T(\xi) v$ by $A^{p-1}$, with $p \in \R$, we obtain the equality
\begin{equation*}
    \frac{1}{h} a(\xi) A^{p-1} v = - c(\xi) A^p v + b(\xi) A^{p-1} B v.
\end{equation*}
Therefore, denoting $\left < z_1, z_2 \right > := z_1^* z_2$ for all $z_1,z_2 \in \C^N$ and setting
\begin{equation}
    \label{eq:3:ecuacionDefinicionDeMu,y}
    y = -h \frac{\left < v , A^p v \right >}{\left < v , A^{p-1} v \right >}, \quad
    \mu = \frac{\left < v , A^{p-1} B v \right >}{\left < v , A^p v \right >},
\end{equation}
we have that
\begin{equation}
    a(\xi) = y c(\xi) - y \mu b(\xi).
    \label{eq:3:ecuacionEstabilidad}
\end{equation}
In addition, since $A$ is positive definite, it is verified that $y<0$ for all $v$, which gives meaning to the following definitions:

\begin{definition}
    For a given $y \in \R^-$ and $\mu \in \C$, the equation (\ref{eq:3:ecuacionEstabilidad}) is stable if every solution satisfies the condition $|\xi| < 1$.
\end{definition}

\begin{definition}
    The equation (\ref{eq:3:ecuacionEstabilidad}) is stable for $y = -\infty$ and $\mu \in \C$ if every solution $\xi$ of the equation $c(\xi)- \mu b(\xi) = 0$ verifies that $|\xi| < 1$.
\end{definition}

\begin{definition}
    We define the region of unconditional stability $D$ as the values of $\mu$ such that the polynomial equation (\ref{eq:3:ecuacionEstabilidad}) is stable for all $y \in \R^- \cup \{ -\infty \}$; that is,
    \begin{equation*}
        D = \bigcap \limits_{y \in \R^- \cup \{ -\infty \}} D_y,
    \end{equation*}
    with $D_y := \{ \mu \in \C :$ (\ref{eq:3:ecuacionEstabilidad}) is stable for $y\}$, $y \in \R^- \cup \{ -\infty \}.$
\end{definition}

\begin{proposition}
    \label{prop:3:caracterizaciónRegiónEstabilidadIncondicional}
    If $u=0$ and $\theta > 1/2$, the region of unconditional stability is
    \begin{equation}
        D = D_{-\infty} = \{ \mu \in \C: |\mu| <1 \}.
        \label{eq:3:ecuacionRegionEstabilidadIncondicional}
    \end{equation}
\end{proposition}

\begin{proof}
    We first show that $D_{-\infty} = D(0,1) := \{ \mu \in \C: |\mu| <1 \}$. It is verified that $c(z)- \mu b(z) = 0$ if and only if $\mu = z^m$ or $z = 1- \frac{1}{\theta}$. Then, since $\big|1-\frac{1}{\theta}\big|<1$ for all $\theta > 1/2$, it is clear that
    \begin{equation*}
        \mu \in D_{-\infty} \Longleftrightarrow \mu \neq z^m \text{ for all } |z| \geq 1 \Longleftrightarrow \mu \in D(0,1).
    \end{equation*}

    Hence, to prove (\ref{eq:3:ecuacionRegionEstabilidadIncondicional}), it suffices to show that $D_{-\infty} \subseteq D_y$ for all $y < 0$. We denote
    \begin{equation*}
        \varGamma_y = \left \{ \frac{1}{b(z)} \left ( c(z) - \frac{a(z)}{y}\right ) : |z| = 1\right \}
    \end{equation*}
    and define
    \begin{equation*}
        P(z, \Tilde{y}) := c(z) - \mu b(z) + \frac{\Tilde{y}}{1-\Tilde{y}} a(z).
    \end{equation*}
    Therefore, for a given $y \in \R^- \cup \{ -\infty \}$, it is verified that $D_y = \left \{ \mu \in \C : P \left (z, \frac{1}{1-y} \right ) = 0 \Longrightarrow |z| < 1 \right \}$ and, by the argument principle,
    \begin{equation*}
        \mu \in D_y \Longleftrightarrow P \left (z, \frac{1}{1-y} \right) \text{ have } m+1 \text{ roots in } D(0,1)
        \Longleftrightarrow
        \begin{cases}
            P \left (z, \frac{1}{1-y} \right) \neq 0 \quad \forall |z| = 1, \\ m+1 = F \left (\frac{1}{1-y} \right),
        \end{cases}
    \end{equation*}
    where
    \begin{equation*}
        F(\Tilde{y}) := \frac{1}{2 \pi i} \int_{|z| = 1} \frac{P_z(z,\Tilde{y})}{P(z,\Tilde{y})} dz, \quad P_z(z,\Tilde{y}) := \frac{\partial P(z,\Tilde{y})}{\partial z}.
    \end{equation*}

    Now, we suppose, by reduction to the absurd, that $\mu \in D_{-\infty} $, but $\mu \notin D_{y_0} $, $y_0 \in \R^-$. Since $\frac{1}{1-y} < 1 $ for all $y < 0$ and $P(z,\Tilde{y})$ is a continuous function for all $\Tilde{y}<1$, then it exists a value $\Tilde{y}_1 \in \left (0, \frac{1}{1-y_0} \right ]$ such that $P(z, \Tilde{y}_1) = 0$ for some $|z| = 1$. That is, $\mu \in \varGamma_{y_1}$ for some $y_1 \in \R^-$ such that $y_1<y_0$.

    However, for a given $y \in \R^-$, we have that
    \begin{equation}
        \begin{split}
            \Gamma_y &= \{ \mu \in \C : a(e^{i \alpha})= yc(e^{i \alpha}) - y \mu b(e^{i \alpha}), \: -\pi \leq \alpha \leq \pi \}\\
            &= \left \{ \frac{ \left(e^{i\alpha} \right)^{m}(1 - e^{i \alpha} + y (1-\theta + e^{i \alpha}\theta))}{y (1 - \theta + e^{i \alpha } \theta )} , \: -\pi \leq \alpha \leq \pi \right \}.
        \end{split}
        \label{eq:3:caracterizaciónGamma_y}
    \end{equation}
    Therefore, if $\mu \in \varGamma_{y} $, it is verified that
    \begin{equation*}
        |\mu| = \left |1+\frac{ 1 - e^{i \alpha}}{y (1 - \theta + e^{i \alpha } \theta)} \right| \geq 1.
    \end{equation*}
    That is, $D(0,1) \cap \varGamma_{y} = \emptyset$ for all $y \in \R^-$, with what is concluded.
    \qed
\end{proof}

\begin{remark}
    If $u=0$ and $0 \leq \theta \leq 1/2$, then the region of unconditional stability is $D = \emptyset$.
\end{remark}

Now, we consider the set $F_p := \left \{ \left < v , A^{p-1} B v \right > : \left < v , A^p v \right > = 1 \right \}$. Then, since the matrix $A$ is positive definite, it can be proved that
\begin{equation*}
    \mu = \frac{\left < v , A^{p-1} B v \right >}{\left < v , A^p v \right >} \in F_p,
\end{equation*}
and, using a change of variables of the type $v = A^{-\frac{p}{2}} x$, we can write
\begin{equation*}
    F_p = F \left (A^{\frac{p}{2}-1} B  A^{-\frac{p}{2}} \right) = \left \{ x^*  A^{\frac{p}{2}-1} B  A^{-\frac{p}{2}}x : x^*x = 1, \: x \in \C^N \right \}.
\end{equation*}
Therefore, by Proposition \ref{prop:3:caracterizaciónRegiónEstabilidadIncondicional}, we can guarantee the stability of (\ref{eq:3:thetaMetodoEstabilidad1}) under certain conditions given by the following theorem:

\begin{theorem}
    \label{thm:3:EstabilidadThetaMetodo}
    If $u=0$, $\theta > 1/2$ and it exists a point $p \in \R$ such that $F \left (A^{\frac{p}{2}-1} B  A^{-\frac{p}{2}} \right) \subseteq D(0,1)$, the $\theta$-method (\ref{eq:3:thetaMetodoEstabilidad1}) is unconditionally stable.
\end{theorem}

\begin{remark}
    If the conditions of the previous theorem are verified, then the corresponding $\theta$-method will be stable regardless of the step size $h = \tau /m$, $m \in \N$.
\end{remark}

\begin{remark}
     Theorem \ref{thm:3:EstabilidadThetaMetodo} cannot be applied if there is no point $p \in \R$ such that $\sigma \left(A^{\frac{p}{2}-1} B  A^{-\frac{p}{2}} \right) \subseteq D(0,1)$.
\end{remark}

On the other hand, it is also verified that if it exists a value $p \in \R$ such that $F \left (A^{\frac{p}{2}-1} B  A^{-\frac{p}{2}} \right) \subseteq  D_y$, where $y$ is as before, then (\ref{eq:3:thetaMetodoEstabilidad1}) is stable. Nevertheless, $y$ depends on a parameter $v$ that, in general, we don't know.
For this reason, to apply the above result, we must consider the inclusion $y \in -h F(A)$, which is clear if we set $v = A^{-\frac{p-1}{2}}x$.
In this way, we have the following proposition, which presents a weaker condition for the stability of (\ref{eq:3:thetaMetodoEstabilidad1}) than Theorem \ref{thm:3:EstabilidadThetaMetodo}:

\begin{proposition}
    \label{prop:3:EstabilidadThetaMetodo}
    If it exists a point $p \in \R$ such that $F \left (A^{\frac{p}{2}-1} B  A^{-\frac{p}{2}} \right) \subseteq \bigcap \limits_{y \in  -hF(A)} D_y$, the $\theta$-method (\ref{eq:3:thetaMetodoEstabilidad1}) is stable.
\end{proposition}

\begin{remark}
    Unlike Theorem \ref{thm:3:EstabilidadThetaMetodo}, Proposition \ref{prop:3:EstabilidadThetaMetodo} is not limited to the study of the stability of the $\theta$-methods with $u=0$ and $\theta > 1/2$. However, this result will indeed depend on the step size considered.
\end{remark}

\subsection{Characterization of $D_y$ }
\label{sec:3.1}

As we have already seen, to study the stability of the $\theta$-methods for the DDE (\ref{eq:1:EDR_EstabilidadThetaMétodo}), we have to check the inclusion $\mu \in D_y$, where $\mu$ and $y$ are the parameters given by (\ref{eq:3:ecuacionDefinicionDeMu,y}). However, to carry out this study using Theorem \ref{thm:3:EstabilidadThetaMetodo} is not always possible.

In this subsection we will provide some results that will facilitate the study of the stability of the $\theta$-method with $u=0$ and $\theta = 1$ when the conditions of Theorem \ref{thm:3:EstabilidadThetaMetodo} are not satisfied. We will simplify the computation of the region $D_y$ and analyze the behavior of this region according to $y$, proving that if $y_1 < y_2 < 0$, then $D_{y_1} \subseteq D_{y_2} $.
In what follows, we will assume the values $u=0$ and $\theta = 1$ to be fixed.

\begin{figure}[htb]
    \centering
    \begin{minipage}[t]{0.47\textwidth}
        \centering
        \includegraphics[width=\textwidth,height=6cm]{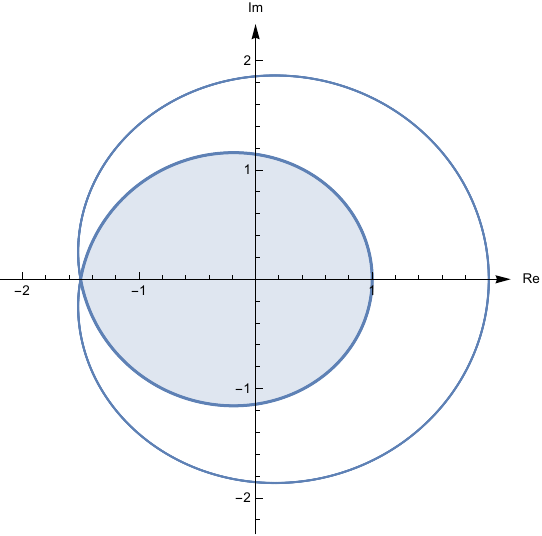}
        (a) Case $m=2$\phantom{xxxx}
    \end{minipage}
    \hfill
    \begin{minipage}[t]{0.47\textwidth}
        \centering
        \includegraphics[width=\textwidth,height=6cm]{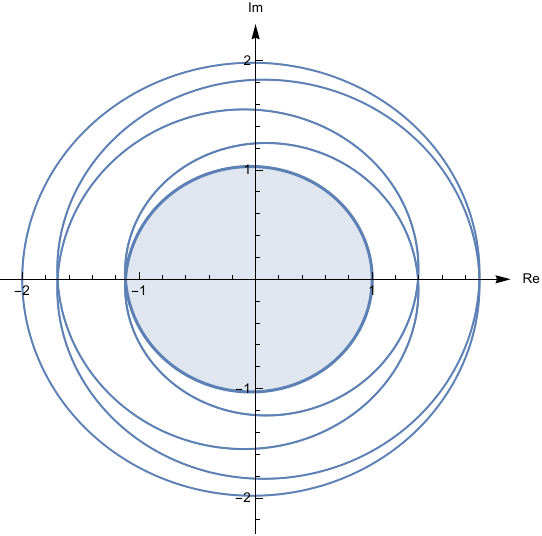}
        (b) Case $m=5$\phantom{xxxx}
    \end{minipage}
    \caption{$\Gamma_y$ for $m=2$ and $m=5$ if $u=0$, $\theta = 1$ and $y = -2$}
    \label{fig:3.1:Gamma_y}
\end{figure}

Let $\Tilde{D}_y$ be the innermost region bounded by the curve $\Gamma_y $. Then, for all $y<0$ it is verified that $\Tilde{D}_y$ is the same region as $D_y $. By way of example, taking into account that
Fig. \ref{fig:3.1:Gamma_y} shows, in the complex plane, the curve $\Gamma_y$ for $y = -2$ and different values of $m$, we are asserting that $D_y$ is, in each case, the region shaded therein\footnote{All figures in this paper have been created using codes written in Mathematica, specifically in Wolfram Mathematica version 13.3.}.

To prove the above statement we first have to see that the region $\Tilde{D}_y$ always exists. For this reason, we are going to demonstrate that $\Gamma_y$ is a closed curve.

Given a value $y \in \R^-$, we can write the curve $\Gamma_y$ with $u=0$ and $\theta=1$ as
\begin{equation*}
    \Gamma_y = \left \{ \frac{1}{y} \left (e^{i\alpha} \right)^{m - 1}\left (1 + e^{i\alpha} (y - 1) \right), \: -\pi \leq \alpha \leq \pi \right \}.
\end{equation*}
Then, the real part of any $\mu \in \Gamma_y$ is given by the expression
\begin{equation}
    \Re(\mu) = \frac{1}{y}\left[\cos ((m-1)\alpha)+(y-1) \cos (m\alpha)\right],
    \label{eq:3.1:parteRealFronteraRegionInterna}
\end{equation}
whereas the imaginary part is
\begin{equation}
    \Im(\mu) = \frac{1}{y} \left [\sin ((m-1)\alpha)+(y-1) \sin (m \alpha) \right].
    \label{eq:3.1:parteImaginariaFronteraRegionInterna}
\end{equation}

Therefore, by the odd symmetry of the sine function and the even symmetry of the cosine, the curve $\Gamma_y$ is symmetric about the real axis.
But moreover, both setting $\alpha = 0$ and $\alpha = \pi$ we have that $\Im(\mu) = 0$. Hence, it is clear that $\Gamma_y$ is a continuous and closed curve, as we wanted to prove.\\

On the other hand, if we denote by $\mu(\alpha,y)$ the point $\mu \in \Gamma_y$ given by $\xi = e^{i \alpha}$, $\alpha \in [-\pi,\pi]$,  we have the following proposition:

\begin{proposition}
    \label{prop:3.1:moduloGamma_yParaDistintosAlpha}
    Let $\alpha_1,\alpha_2 \in [ -\pi, \pi]$. If $y \in \R^-$ and $|\alpha_1| < |\alpha_2|$, it is verified that $|\mu(\alpha_1,y)|< |\mu(\alpha_2,y)|$.
\end{proposition}

\begin{proof}
    Since
    \begin{equation}
        \mu(\alpha,y)=\frac{1}{y} \left (e^{i\alpha} \right)^{m - 1}\left (1 + e^{i\alpha} (y - 1) \right)= \left (e^{i\alpha} \right)^{m - 1} \left (1 + \left(1-\frac{1}{y}\right) (-1+e^{i\alpha})  \right )
        \label{eq:3.1:caracterizacionMu(alpha,y)}
    \end{equation}
    and $\big |e^{i \alpha} \big|=1$ for all $\alpha \in [-\pi, \pi]$, we only have to compare the modules of the expressions $1+\phi_1$ and $1+\phi_2 $, where
    $\phi_1 :=  \left(1-\frac{1}{y}\right) (-1+e^{i\alpha_1})$ and $\phi_2 := \left(1-\frac{1}{y}\right) (-1+e^{i\alpha_2}) $.

    Therefore, since for all $z \in \C$ it is verified that
    \begin{equation}
        |1+z|^2 = 1 + |z|^2 + 2 \Re(z),
        \label{eq:3.1:propiedadModulo}
    \end{equation}
    it will suffice to prove that $|\phi_1|^2 + 2 \Re(\phi_1) < |\phi_2|^2 + 2 \Re(\phi_2)$.

    However, we can write these expressions as
    \begin{equation*}
        \begin{split}
            |\phi_1|^2 + 2 \Re(\phi_1) &= \left(1-\frac{1}{y}\right)^2 \left [(-1+\cos(\alpha_1))^2+\sin(\alpha_1)^2 \right ] + 2\left(1-\frac{1}{y}\right) (-1+ \cos(\alpha_1))\\
            &= 2 \left(1-\frac{1}{y}\right)^2 (1-\cos(\alpha_1))-2 \left(1-\frac{1}{y}\right)(1-\cos(\alpha_1))\\
            &=2 \left(1-\frac{1}{y}\right)\left(-\frac{1}{y}\right)(1-\cos(\alpha_1))
        \end{split}
    \end{equation*}
    and
    \begin{equation*}
        |\phi_2|^2 + 2 \Re(\phi_2) = 2 \left(1-\frac{1}{y}\right)\left(-\frac{1}{y}\right)(1-\cos(\alpha_2)).
    \end{equation*}
    Hence, since $\cos(\alpha_1) > \cos(\alpha_2)$, it concludes.
    \qed
\end{proof}

Having seen the above result, we are now ready to test the inclusion $\Tilde{D}_y \subseteq D_y$.
For it, we will study the behavior of the real part and the imaginary part of $\varGamma_y$ as a function of $\alpha$.

Since $\mu(0,y) =1$, we start by analysing the behavior of $\Im(\mu(\alpha,y))$ in the half-plane $\{ z \in \C : \Re(z) \geq 0 \}$. It follows from
(\ref{eq:3.1:parteRealFronteraRegionInterna}) that if $\Re(\mu(\alpha,y)) \geq 0$, the derivative of $\Im(\mu(\alpha,y))$ with respect to $\alpha$ verifies
\begin{equation*}
       \frac{\partial}{\partial \alpha} \Im(\mu(\alpha,y)) = \frac{1}{y}[(m-1) \cos ((m-1)\alpha  )+m (y-1) \cos(m \alpha )]  \geq \left(1-\frac{1}{y}\right) \cos(m \alpha).
\end{equation*}
Therefore, for any sufficiently small $\alpha >0$, we will have that $\Im(\mu(\alpha,y)) > 0$.

Suppose then that $\Im(\mu(\alpha,y)) > 0$. In that case, we can deduce easily from (\ref{eq:3.1:parteImaginariaFronteraRegionInterna}) that
\begin{equation*}
       \frac{\partial}{\partial \alpha} \Re(\mu(\alpha,y)) = -\frac{1}{y}[(m-1) \sin ((m-1)\alpha  ) + m (y-1) \sin(m \alpha )] < \frac{1}{y} \sin((m-1) \alpha).
\end{equation*}
That is, as long as $\Im(\mu(\alpha,y)) > 0$ and $\alpha \in \big [0, \frac{\pi}{m-1} \big ]$, the real part $\Re(\mu(\alpha,y))$ will be decreasing with respect to $\alpha$.
Hence, by Proposition \ref{prop:3.1:moduloGamma_yParaDistintosAlpha} and since  $\Im \big( \mu \big (\frac{\pi}{m-1},y \big ) \big )\leq 0$ for all $m \in \N-\{1\}$, we can conclude that $D(0,1) \subseteq \Tilde{D}_y$.

But moreover, it is clear that the roots of (\ref{eq:3:ecuacionEstabilidad}) when $\mu = 0$ and $\theta=1$ are $\xi = 0$ (with multiplicity $m$) and $\xi = \frac{1}{1-y}$. Therefore, $0 \in D_y$ and, by continuity, $\Tilde{D}_y \subseteq D_y$. \\

We will now prove the opposite inclusion.
For each $y \in \R^-$ we denote by $\gamma_y$ the curve $\gamma_y : [-\pi, \pi] \rightarrow \varGamma_y \subseteq \C$ given by $\gamma_y(\alpha) = \mu(\alpha,y)$. Then, by studying its argument, it can be concluded that $\gamma_y$ turns counterclockwise around the point 0 during its entire domain.
But furthermore, it can be observed that the index of $0$ with respect to $\gamma_y$ is
\begin{equation*}
    W_{\gamma_y}(0)= \frac{1}{2 \pi i} \int_ {\gamma_y} \frac{1}{z} dz =  m.
\end{equation*}
Therefore, the curve $\gamma_y$ turns exactly $m$ counterclockwise around the point $0$.

On the other hand, by Proposition \ref{prop:3.1:moduloGamma_yParaDistintosAlpha}, we know that if we have two values $\alpha_1, \alpha_2 \in [-\pi, \pi]$ such that $\Im(\mu(\alpha_1,y)) \neq 0$ and $\alpha_1 \neq \alpha_2$, then $\mu(\alpha_1,y) \neq \mu(\alpha_2,y)$.
But moreover, if $\Re(\mu(\alpha,y))=0$, it can be proved that the root $e^{i \alpha}$ of the equation
\begin{equation*}
    P(\alpha, y, \xi) := a(\xi) - y c(\xi) + y \mu(\alpha,y) b(\xi) = 0
\end{equation*}
has multiplicity $1$, since for every $\mu(\alpha,y)$ with zero real part it is verified that
\begin{equation*}
    \Re \left (\frac{\partial}{\partial \xi} P(\alpha, y, e^{i \alpha}) \right ) = 0 \Longleftrightarrow \cos ((m-1)\alpha) = 0 \Longleftrightarrow \cos (m \alpha) = 0,
\end{equation*}
and, however, the imaginary part of the above derivative is not null whenever $\cos ((m-1)\alpha) = \cos (m \alpha) = 0$.
Therefore, for each point $\mu(\alpha,y)$ with zero real part we have that the equation (\ref{eq:3:ecuacionEstabilidad}) has no root with multiplicity longer than $1$.

But furthermore, for all fixed $y<0$ it follows that if $\mu$ is sufficiently large, then (\ref{eq:3:ecuacionEstabilidad}) has $m$ roots with modulus higher than $1$.
Hence, we can conclude that $D_y \subseteq \Tilde{D}_y$ and, consequently, that $D_y = \Tilde{D}_y$. \\

We now present other results that we will need to prove Theorem \ref{thm:3.1:comparaciónDy}, which compares the regions $D_y$ given by two different parameters $y$.

\begin{remark}
    Since $\Im \left (\mu \left (\frac{\pi}{m},y \right ) \right) \leq 0$ for all $y \in \R^-$, we can write the boundary of $D_y$ as $$\partial D_{y}=\{ \mu(\alpha,y) \in \Gamma_{y} : \alpha \in [-\beta,\beta] \},$$
    where $\beta$ is a certain value belonging to the interval $\left[0,\frac{\pi}{m} \right]$. Hence, given an arbitrary point $\mu(\alpha,y) \in \partial D_y$, we have:
    \begin{equation}
        \begin{split}
            &\Im(\mu(\alpha,y))>0 \Longleftrightarrow \alpha \in (0,\beta)
            , \\ &\Im(\mu(\alpha,y))<0 \Longleftrightarrow \alpha \in (-\beta,0).
        \end{split}
        \label{eq:3.1:relacionIm(mu(alpha,y))ConAlpha}
    \end{equation}
\end{remark}

\begin{proposition}
    \label{prop:3.1:comparacionMu(alpha,y1)Mu(alpha,y2)}
    Let $y_1,y_2 \in \R^-$. If $y_1<y_2<0$ , then $|\mu(\alpha,y_1)|< |\mu(\alpha,y_2)|$ for all $\alpha \neq 0$.
\end{proposition}

\begin{proof}
    By (\ref{eq:3.1:caracterizacionMu(alpha,y)}),
    we only have to prove that $|1 + \varphi_1| < |1 + \varphi_2|$, where $\varphi_1 :=  \left(1-\frac{1}{y_1}\right) (-1+e^{i\alpha})$ and $\varphi_2 := \left(1-\frac{1}{y_2}\right) (-1+e^{i\alpha})$.
    Nevertheless, if $\alpha \neq 0$, then
    \begin{equation*}
        \begin{split}
            |\varphi_1|^2 + 2 \Re(\varphi_1) &= 2 \left(1-\frac{1}{y_1}\right)\left(-\frac{1}{y_1}\right)(1-\cos(\alpha)) \\ &< 2 \left(1-\frac{1}{y_2}\right)\left(-\frac{1}{y_2}\right)(1-\cos(\alpha)) = |\varphi_2|^2 + 2 \Re(\varphi_2).
        \end{split}
    \end{equation*}
    Therefore, by (\ref{eq:3.1:propiedadModulo}), it concludes.
    \qed
\end{proof}

Taking all these results into account, we can already prove the following theorem:

\begin{theorem}
    \label{thm:3.1:comparaciónDy}
    If $y_1 < y_2 <0$, then $D_{y_1} \subseteq D_{y_2}$.
\end{theorem}

\begin{proof}
    We assume first that $m=1$. In that case, $\mu(\alpha,y) = \frac{1}{y} + e ^{i \alpha} \left(1-\frac{1}{y} \right)$ and, consequently, the curve $\partial D_y = \Gamma_y$, which coincides with the circumference with center $\frac{1}{y}$ and radius $1-\frac{1}{y}$. Therefore, since $1 \in \partial D_{y_1}, \partial D_{y_2}$, if $m=1$, it suffices to compare the center and the radius of both boundaries to conclude that $D_{y_1} \subseteq D_{y_2}$.

    Therefore, if, under the previous hypothesis, $\mu(\alpha_1,y_1) \in \partial D_{y_1}$ and $\mu(\alpha_2,y_2) \in \partial D_{y_2}$ are two points with the same argument, then $|\mu(\alpha_1,y_1)| \leq |\mu(\alpha_2,y_2)|$. That is, for a point $\mu(\bar \alpha_1,y_1) \in \partial D_{y_1}$ to have the same module as $\mu(\alpha_2,y_2) \in \partial D_{y_2}$ when $m=1$, the constraint $|\bar \alpha_1|\geq |\alpha_1|$ must be verified and, consequently, if $\alpha_1, \alpha_2 >0$, the condition $\arg(\mu(\bar \alpha_1,y_1)) \geq \arg(\mu(\alpha_1,y_1)) = \arg(\mu(\alpha_2,y_2))$, where $\arg: \C-{0} \rightarrow [0, 2 \pi)$, must also be fulfilled.

    From the previous reasoning, we will see now that the inclusion of the statement is also verified for $m \in \N - \{ 1 \}$. We assume again that $\mu(\alpha_1,y_1) \in \partial D_{y_1}$ and $\mu(\alpha_2,y_2) \in \partial D_{y_2}$ are two points with the same argument, but now setting $m \neq 1$.
    We have to see that if $\alpha_1$ and $\alpha_2$ are not null, then $\mu(\alpha_1,y_1) \neq \mu(\alpha_2,y_2)$; that is, $|\mu(\alpha_1,y_1)| \neq |\mu(\alpha_2,y_2)|$.

    By (\ref{eq:3.1:relacionIm(mu(alpha,y))ConAlpha}), we can assume that $\alpha_1$ and $\alpha_2$ share the same sign. Therefore, by symmetry, it will suffice to show that $|\mu(\alpha_1,y_1)| \neq |\mu(\alpha_2,y_2)|$ for $\alpha_1, \alpha_2 \in \left(0 , \frac{\pi}{m} \right]$.
    We suppose, by reduction to the absurd, that $|\mu(\alpha_1,y_1)| = |\mu(\alpha_2,y_2)|$. Then, by Proposition \ref{prop:3.1:moduloGamma_yParaDistintosAlpha} and Proposition \ref{prop:3.1:comparacionMu(alpha,y1)Mu(alpha,y2)}, whenever $\alpha_1$ and $\alpha_2 $ are higher than $0$, it must the condition $\alpha_1 > \alpha_2$ is met.

    But on the other hand, given a value $y \in \R^-$, we have that
    \begin{equation}
        \arg(\mu(\alpha,y)) = \arg \left (e ^{i (m-1)\alpha} \right ) + \arg \left(\frac{1}{y} + e ^{i \alpha} \left(1-\frac{1}{y} \right)\right) +2 k \pi, \quad k \in \Z,
        \label{eq:3.1:caracterizaciónArgumentoMu(alpha,y)}
    \end{equation}
    and, moreover, given $\alpha \in \left(0, \frac{\pi}{m} \right]$, it can also be verified that the arguments $ \arg \left (e ^{i (m-1)\alpha} \right )$ and $\arg \left(\frac{1}{y} + e ^{i \alpha} \left(1-\frac{1}{y} \right)\right)$ belong to the interval $(0,\pi)$ whenever $m \in \N - \{ 1 \}$. Hence, for $\alpha = \alpha_1, \alpha_2$ we can write
    \begin{equation*}
         \arg (\mu(\alpha,y)) = (m-1)\alpha + s(\alpha,y), \quad s(\alpha,y) := \arg \left(\frac{1}{y} + e ^{i \alpha} \left(1-\frac{1}{y} \right)\right),
    \end{equation*}
    which implies that $s(\alpha_1,y_1) \geq s(\alpha_2,y_2)$, since the module of $\mu(\alpha,y)$ does not depend on the value $m$ and the argument $\arg(\mu(\alpha,y)) = s(\alpha,y)$ when $m=1$. That is, we have that $\arg(\mu(\alpha_1,y_1)) > \arg(\mu(\alpha_2,y_2))$, which contradicts our starting hypothesis.

    Hence, $\partial D_{y_1} \cap \partial D_{y_2} = {1}$ and, as we know that $D_{-\infty} \subseteq D_y$ for all $y \in \R^-$ , by continuity the theorem is proved.
    \qed
\end{proof}

Combining the previous theorem, and Proposition \ref{prop:3:EstabilidadThetaMetodo}, we obtain the following result:

\begin{theorem}
    \label{prop:3:EstabilidadThetaMetodo2}
    Let us consider that $\lambda_j$ is the highest eigenvalue of $A$, and $   y_j = -h \lambda_j$.  If it exists a point $p \in \R$ such that $F \left (A^{\frac{p}{2}-1} B  A^{-\frac{p}{2}} \right) \subseteq   D_{y_j}$, the $\theta$-method (\ref{eq:3:thetaMetodoEstabilidad1}) is stable.
\end{theorem}

\subsection{ Relation with other works}
\label{sec:3.2}

In other studies \cite{Calvo,Rihan}, the stability of $\theta$-methods was analyzed employing Eq. (\ref{eq:1:EDR_EstabilidadTipica}), and assuming   that the matrices $A$ and $B$, in DDE (\ref{eq:1:EDR_EstabilidadThetaMétodo}), are simultaneously diagonalizable.  This idea is well-known, and it was also employed in a large number of papers for studying the stability of numerical methods for systems of ODEs and PDEs without delay, as it was mentioned in the introduction.

Let us consider the test equation
 \begin{equation}
   \label{eq:1:EDR_EstabilidadTipica2}
   y'(t) = \lambda ( y(t) + \mu y(t-\tau)  ),
\end{equation}
where $\lambda, \gamma \in \C$, $\mu = \frac{\gamma}{\lambda}$.  From the standard analysis, the following result is straightforward  when matrices $A$ and $B$ are simultaneously diagonalizable:

\begin{proposition}
    For every $i=1,\ldots,N$, where $N$ is the dimension of the matrices $A$ and $B$, let $\lambda_i$ and $\gamma_i$ be the eigenvalues of $A$ and $B$, respectively, corresponding to the same eigenvector $v_i$. Then, a certain $\theta$-method is stable (for a step size $h>0$) whenever for every $i$ we have that $\mu_i := \frac{\gamma_i}{\lambda_i} \in D_{y_i}$, $y_i := -\lambda_i h$.
\end{proposition}

\begin{remark}
    By the previous proposition, if it exists some $i$ such that $\Re(\mu_i) \geq 1$, the $\theta$-method with $u=0$ and $\theta = 1$ is unconditionally unstable; that is, it is not stable for any step size.
\end{remark}

Now, we will use the results obtained earlier to simplify the study of stability when matrices $A$ and $B$ are simultaneously diagonalizable, therefore the matrix $A$ will need not be symmetrical, but it will be sufficient that all its eigenvalues are all real and higher than $0$. From Theorem \ref{thm:3.1:comparaciónDy}, we can directly obtain:

\begin{proposition}
\label{prop:4:condiciónSuficienteEstabilidadCondicional}
    Let $\lambda_i$, $\mu_i$ and $y_i$ be as before and assume that $\lambda_j$ is the highest eigenvalue of $A$. If for every $i$ we have $\mu_i \in D_{y_j}$, then the $\theta$-method with $u=0$ and $\theta=1$ is stable.
\end{proposition}

\begin{remark}
   When matrices $A$ and $B$ are simultaneously diagonalizable in $\R$, the result above is equivalent to Proposition \ref{prop:3:EstabilidadThetaMetodo2}.  Since both matrices are  simultaneously diagonalizable,   $F \left (A^{-1} B \right) $ becomes the interval $[\mu_{1},  \mu_{2}]$, being $\mu_1$ the smallest generalized eigenvalue, and $\mu_{2}$ the largest one.  And, obviously, Proposition \ref{prop:3:condiciónSuficienteEstabilidadConMatricesSimultáeamenteDiagonalizables1} (which comes directly from Proposition \ref{prop:4:condiciónSuficienteEstabilidadCondicional}) is equivalent to Theorem \ref{thm:3:EstabilidadThetaMetodo}.
\end{remark}

\begin{proposition}
    \label{prop:3:condiciónSuficienteEstabilidadConMatricesSimultáeamenteDiagonalizables1}
    If the condition $|\mu_i| < 1$ is verified for each $i$, where $\mu_i$ is as before, then all $\theta$-method with $u=0$ and $\theta>1/2 $ is unconditionally stable (regardless of the step size considered).
\end{proposition}

In \cite[Th.~3.1]{kj1994stability}, K.~in't Hout studied the unconditional stability of $\theta$-methods in a more general case ($A$ and $B$ do not need to be simultaneously diagonalizable).  Let us define
$H^*:= \{ \xi: \xi \in \C, Re( \xi) < 0\}$, he demonstrated:

\begin{theorem}
    \label{prop:4:Karel}
    Let us consider: \\
    (i)  $\sigma(A) \subset H^*$,  $\sup_{ Re(\xi)=0 } \rho\left( (\xi I -A)^{-1} B \right) < 1$,
    and $-1 \notin  \sigma(A^{-1} B )$, \\
   (ii)    The $\theta$-method (\ref{eq:3:thetaMetodoEstabilidad1}) is unconditionally stable.\\
   (i) and (ii) are equivalent.
\end{theorem}

This theorem above provides a sufficient and enough constrain to obtain unconditional stability, whilst Theorem \ref{thm:3:EstabilidadThetaMetodo} provides only a sufficient condition.  However, if $A$ and $B$ are simultaneously diagonalizable, it is easy to check  that condition (i) in Theorem \ref{prop:4:Karel} is equivalent to the condition proposed in Theorem \ref{thm:3:EstabilidadThetaMetodo}.

Actually, whenever the matrix $A^{-1} B $  is normal,  then $r( A^{-1} B ) = \rho( A^{-1} B )$  \cite[p.~45]{HornJohnson}, where $r(.)$ is the numerical radius defined as follows:
\begin{equation*}
    r(M) := \max \{ |z|: z \in F(M) \}.
\end{equation*}
Additionally,
\begin{equation*}
    \rho( A^{-1} B ) \leq  \sup_{ Re(\xi)=0 } \rho\left( (\xi I -A)^{-1} B \right),
\end{equation*}
hence, whenever condition (i) in Theorem \ref{prop:4:Karel} is satisfied, then   $F \left (A^{-1} B  \right) \subseteq D(0,1)$.  If $\sigma(A ) \subseteq \R^- $, then we have the conditions in Theorem \ref{thm:3:EstabilidadThetaMetodo} to obtain unconditional stability.

On the other hand, if $F \left (A^{\frac{p}{2}-1} B  A^{-\frac{p}{2}} \right) \subseteq D(0,1)$, then the method is unconditionally stable, then condition (i) in Theorem \ref{prop:4:Karel} is satisfied.  Therefore, whenever  $\sigma(A ) \subseteq \R^- $ and $A^{-1} B $  is normal both theorems are equivalent.

However,  in \cite{kj1994stability}, the stability for a given step size is not considered, and to the best of our knowledge, this was not studied, until now, for numerical methods applied to DDEs.  The procedure employed in this paper can be useful for many other more complicated methods, but also different classes of problems.  The topic of the field of values is well-known, and numerous studies analyze its properties.  Also, different algorithms have been developed to calculate it efficiently for large matrices, including the code \textit{fov} in matlab that can be freely download in
\newblock \href{https://www.chebfun.org/}   {\path{https://www.chebfun.org/}}. In this paper, we will use an algorithm similar to the one described in \cite{HornJohnson,Johnson}.
In what follows, we shall explain how the theory developed above can be employed.

\begin{figure}[htb]
    \centering
    \begin{minipage}[t]{0.4\textwidth}
        \centering
        \includegraphics[width=\textwidth,height=6cm]{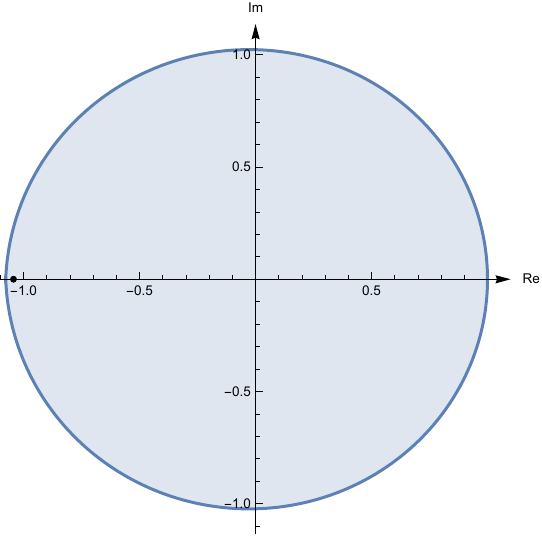}
        (a) Case $m=2$ \phantom{xx}
    \end{minipage}
    \hfill
    \begin{minipage}[t]{0.4\textwidth}
        \centering
        \includegraphics[width=\textwidth,height=6cm]{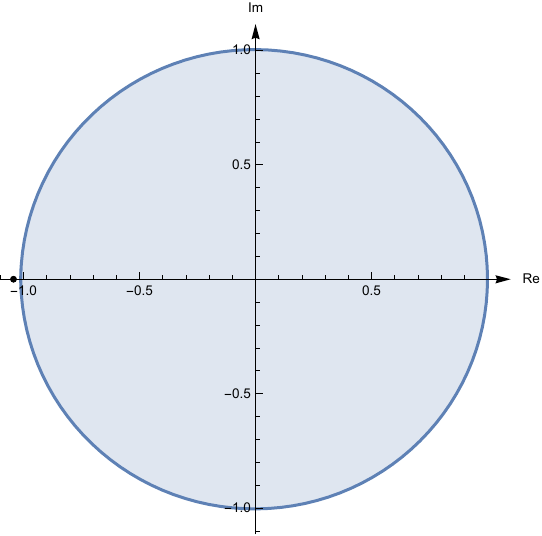}
        (b) Case $m=50$ \phantom{x}
    \end{minipage}
    \caption{Example \ref{ejem:3:MatricesSimultáneamenteDiagonalizables}: comparison of the point $\mu_2$ with the regions $D_{y_1}$ and $D_{y_2}$ given by $m=2$ and $m=50$, respectively}
    \label{fig:3:caracterizacionDyEjemploMatricesSimultáneamenteDiagonalizables}
\end{figure}

\begin{example}
    \label{ejem:3:MatricesSimultáneamenteDiagonalizables}
  Let us study the DDE
    \begin{equation}
        y'(t) = -A \cdot y(t) + B \cdot y(t-\tau)
        \label{eq:3:EDRejemploMatricesSimultáneamenteDiagonalizables}
    \end{equation}
    with
    \begin{equation*}
        A=\begin{pmatrix}
        29 & -7 & 1\\
        3 & 27 & -7\\
        3 & 9 & 11\\
        \end{pmatrix}, \quad
        B=\begin{pmatrix}
        -30 & -27 & 33 \\
        -3 & -96 & 75 \\
        -3 & -111 & 90 \\
        \end{pmatrix}.
    \end{equation*}
    We then have a DDE given by two simultaneously diagonalizable matrices, since the eigenvalues of the matrix $A$, $(\lambda_1, \lambda_2, \lambda_3) = (26,23,18)$, and the eigenvalues of $B $, $(\gamma_1,\gamma_2,\gamma_3) = (-27,-24,15)$, share the same eigenvectors, respectively.
    Nevertheless, in this case we cannot apply Proposition \ref{prop:3:condiciónSuficienteEstabilidadConMatricesSimultáeamenteDiagonalizables1}, so we cannot guarantee the unconditional stability of the $\theta$-methods with $u=0$ and $\theta > 1/2$.

    Even so, these methods can be stable for the DDE (\ref{eq:3:EDRejemploMatricesSimultáneamenteDiagonalizables}) if we consider the appropriate step sizes.
    We suppose the delay is $\tau =1$ and consider the $\theta$-method with $u=0$ and $\theta=1$. In addition, we denote $\mu_i = \gamma_i/\lambda_i$ for every  $i =1,2,3$. It is clear then, by Proposition \ref{prop:3:caracterizaciónRegiónEstabilidadIncondicional}, that for any step size we have $\mu_3 \in D_{y_3}$ (since $|\mu_3| <$1). But moreover, setting $h=0.5$ (that is, $m=2$), we also have $\mu_1,\mu_2 \in D_{y_1}$ and, as a result, the $\theta$-method is stable. However, we can also check that this is not true for any step size. For example, setting $h=0.02$ ($m=50$) we have that $\mu_2 \notin D_{y_2}$ and, therefore, that the $\theta$-method is not stable. Fig. \ref{fig:3:caracterizacionDyEjemploMatricesSimultáneamenteDiagonalizables} shows the point $\mu_2$ together with the region $D_{y_1}$ obtained with $m=2$, first, and together with the region $D_{y_2}$ obtained with $m=50$, later. The inclusion $\mu_1 \in D_{y_1}$ can be easily deduced for $m=2$ from the inequality $\mu_2<\mu_1<0$.

    On the other hand, it can be checked that, for a fixed delay, the solutions of the continuous model (\ref{eq:3:EDRejemploMatricesSimultáneamenteDiagonalizables}) do not have to be asymptotically stable (see \cite{Bellen}).
    However, contrary to what happens with many models and numerical schemes for ODEs and PDEs without delay, in this case the $\theta$-method above considered converges asymptotically for a large step size, while for smaller step sizes it indeed diverges (in line with what is expected from the analytical study for linear systems of DDEs presented in \cite{Bellen}).
\end{example}

\section{Numerical examples}
\label{sec:5}

In this section we show the potential of the theory discussed in this paper by means of two numerical examples given by parabolic partial delay differential equations.

\begin{example}
    Consider first the delayed reaction-diffusion system
    \begin{equation}
        \begin{cases}
            \frac{\partial v_1}{\partial t} = \lambda_1 \frac{\partial^2 v_1}{\partial x^2} - e^{\frac{l\pi}{2}}v_1(t-\tau,x)+ \frac{1}{4} e^{\frac{l\pi}{2}}\left(4 l + \pi^2\right) v_2(t-\tau,x), \quad t \geq 0, \quad 0 \leq x \leq 2,\\
            \frac{\partial v_2}{\partial t} = \lambda_2 \frac{\partial^2 v_2}{\partial x^2} - \frac{1}{4} e^{\frac{l\pi}{2}}\left(4 l + \pi^2 \right)v_1(t-\tau,x)- e^{\frac{l\pi}{2}} v_2(t-\tau,x), \quad t \geq 0, \quad 0 \leq x \leq 2,
        \end{cases}
    \label{eq:4:EDPR1}
    \end{equation}
    with the initial and boundary conditions
    \begin{equation*}
        \begin{cases}
        v_1(t,0)=v_1(t,2)=0,\\
        v_2(t,0)=v_2(t,2)=0,\\
        v_1(t,x) = \phi_1 (t,x),\quad t \leq 0,\\
        v_2(t,x) = \phi_2 (t,x),\quad t \leq 0,
        \end{cases}
    \end{equation*}
    where $\tau$ is a constant delay, $l$ is a real number, $\lambda_1$ and $\lambda_2$ are two positive parameters and $\phi_1$ and $\phi_2$ are two given functions.

    Let $M \in \N$ and $ \Delta x := 2/M$ and define the mesh points $x_0 = 0 < x_1 < \ldots < x_j = j \Delta x < \ldots < x_M=2$. Moreover, for every $i=1,2$ we denote by $z^{i,j}(t)$ the approximation of $v_i(t,x_j)$ and by $z^i(t)$ the vector $(z^{i,1}, \ldots, z^{i,M-1})^T$. Then, by replacing the spatial derivatives with the standard second-order centered differences,
    we obtain the MOL approach
    \begin{equation}
        z'(t) = A z(t) +B z(t-\tau),
        \label{eq:4:EDPR1_Adaptada}
    \end{equation}
    where $z(t)=(z^1(t), z^2(t))^T$,
    \begin{equation*}
        A = \begin{pmatrix}
            \lambda_1 L & 0\\
            0 & \lambda_2 L
        \end{pmatrix}, \quad
        L = \frac{1}{\Delta x^2}
        \begin{pmatrix}
            -2 & 1 & 0 & \ldots & 0\\
            1 & -2 & 1 & \ldots & 0\\
            0 & 1 & -2 & \ldots & 0\\
            \vdots & \vdots & \vdots & \ddots & \vdots\\
            0 & 0 & \ldots & 1 & -2\\
        \end{pmatrix}, \quad
        B= e^{\frac{l \pi}{2}} \begin{pmatrix}
            -I & \left (l + \frac{\pi^2}{4} \right) I\\
            -\left(l+\frac{\pi^2}{4}\right) I & -I\\
        \end{pmatrix}
    \end{equation*}
    and $I$ is the identity matrix of dimension $M-1$.
    That is, we can calculate a numerical solution of (\ref{eq:4:EDPR1}) by applying some $\theta$-method to the equation (\ref{eq:4:EDPR1_Adaptada}).
    We compare the numerical solution obtained in this way with its corresponding exact solution.

    If we set, for example, $\lambda_1 = \lambda_2 = 1$, and $\tau=\frac{\pi}{2}$.
    In that case, it can be verified that the exact solution of the equation (\ref{eq:4:EDPR1}) is
    \begin{equation}
        \begin{cases}
            v_1(t, x) = \phi_1(t,x) =  e^{lt} \sin(t) \sin \left(\frac{\pi x}{2} \right),\\
            v_2(t, x) = \phi_2(t,x) =  e^{lt} \cos(t) \sin \left(\frac{\pi x}{2} \right),
        \end{cases}
        \label{eq:4:solucionExactaEDPR1}
    \end{equation}
    and, as a result, its asymptotic property  changes according to the sign of $l$ (it diverges when $l>0$ and tends to zero when $l<0$).

    To calculate a numerical solution we will use the $\theta$-method with $u=0$, $\theta=1$ and $M=100$ and we will consider the cases $l=-1/10$ and $l=1/10$.
    In the first case, the field of values $F(A^{-1}B)$ is included in the unit disk and, as a result, we can apply Theorem \ref{thm:3:EstabilidadThetaMetodo}; that is, every $\theta$-method with $u=0$ and $\theta>1/2$ is unconditionally stable. Fig. \ref{fig:4:soluciónThetaMétodoEjemplo1} shows the numerical solution given by the above $\theta$-method if we set $l=-1/10$ and $m=5$, in addition to the parameters above defined.
    On the other hand, with $l=1/10$ we can no longer use this Theorem (in fact it can be seen that in that case the $\theta$-method is unstable). Therefore, in both cases we obtain the same results as with the exact solution. Fig. \ref{fig:4:erroresThetaMétodoEjemplo1} shows, together with the unit disk, the field of values $F(A^{-1}B)$ for $l=-1/10$ and $l=1/10$.

    Finally, and again with the above parameters, we present in Table \ref{tab:4:erroresEDPR1} the errors made at $t=10 \pi$ when approximating the exact solution of the problem with the numerical solution obtained for $l=-1/10$ and different values of $m$\footnote{For both $v_1$ and $v_2$ we have calculated this errors from the expression $||e_i||_{2,h}=\sqrt{e_{i,1}^2+ \ldots + e_{i,M-1}^2}$, $i=1,2$, where each $e_{i,j}$ denotes the error made with the corresponding step size $h$ at point $x_j$.}. Note that, as one would expect, as $m$ increases the errors in both components decrease, giving virtually identical solutions for sufficiently large $m$ values.

    \begin{table}[htb]
        \caption{Errors in $t=10 \pi$ of numerical solutions of (\ref{eq:4:EDPR1}) given by $\theta$-method with $u=0$ and $\theta=1$}
        \label{tab:1}
        \begin{tabular}{llllll}
        \hline\noalign{\smallskip}
         &  m=5 & m=25 & m=50 & m=100 & m=1000 \\
        \noalign{\smallskip}\hline\noalign{\smallskip}
        $v_1$ & 0.018354 & 0.006456 & 0.003399 & 0.001697 & 0.000416\\
        $v_2$ & 0.196042 & 0.055879 & 0.029162 & 0.014763 & 0.001122 \\
        \noalign{\smallskip}\hline
        \end{tabular}
        \label{tab:4:erroresEDPR1}
    \end{table}

    \begin{figure}[htb]
        \centering
        \begin{minipage}[t]{0.40\textwidth}
            \centering
            \includegraphics[width=\textwidth,height=5cm]{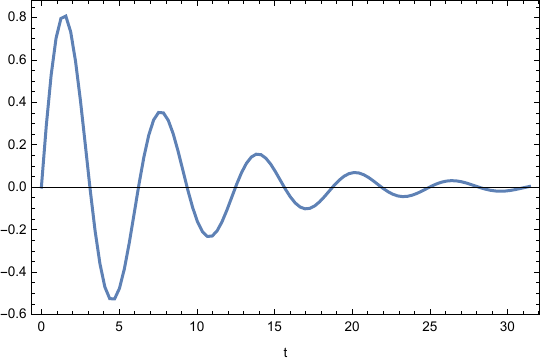}
            \phantom{xx} (a) $v_1$ for $x=1$
        \end{minipage}
        \hfill
        \begin{minipage}[t]{0.55\textwidth}
            \centering
            \includegraphics[width=\textwidth,height=5cm]{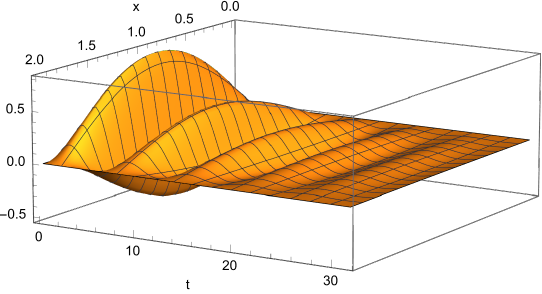}
            (b) $v_1$ along the whole interval
        \end{minipage}

        \vspace{0.5cm}

        \begin{minipage}[t]{0.40\textwidth}
            \centering
            \includegraphics[width=\textwidth,height=5cm]{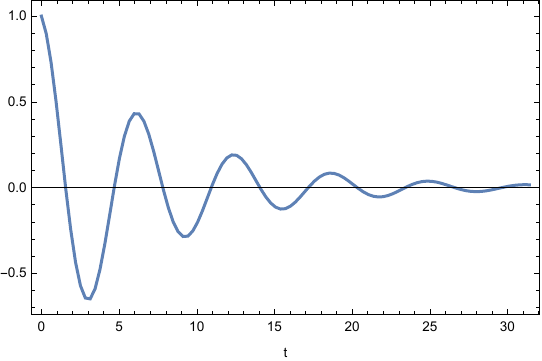}
            \phantom{xx} (c) $v_2$ for $x=1$
        \end{minipage}
        \hfill
        \begin{minipage}[t]{0.55\textwidth}
            \centering
            \includegraphics[width=\textwidth,height=5cm]{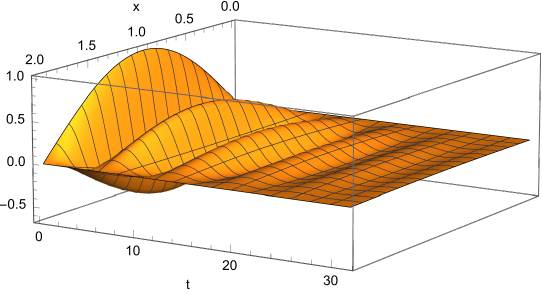}
            (d) $v_2$ along the whole interval
        \end{minipage}

        \caption{Numerical solution of (\ref{eq:4:EDPR1}) given by the $\theta$-method with $u=0$ and $\theta=1$}

        \label{fig:4:soluciónThetaMétodoEjemplo1}
    \end{figure}

    \begin{figure}[htb]
        \centering
        \begin{minipage}[t]{0.4\textwidth}
            \centering
            \includegraphics[width=\textwidth,height=5cm]{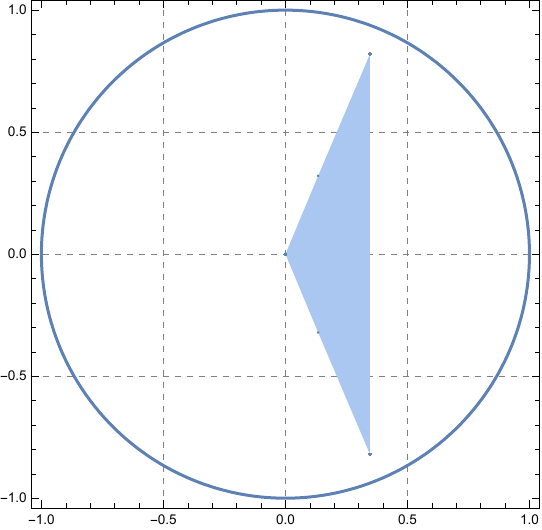}
            \phantom{xx} (a) Case $l=-1/10$
        \end{minipage}
        \hspace{2cm}
        \begin{minipage}[t]{0.4\textwidth}
            \centering
            \includegraphics[width=\textwidth,height=5cm]{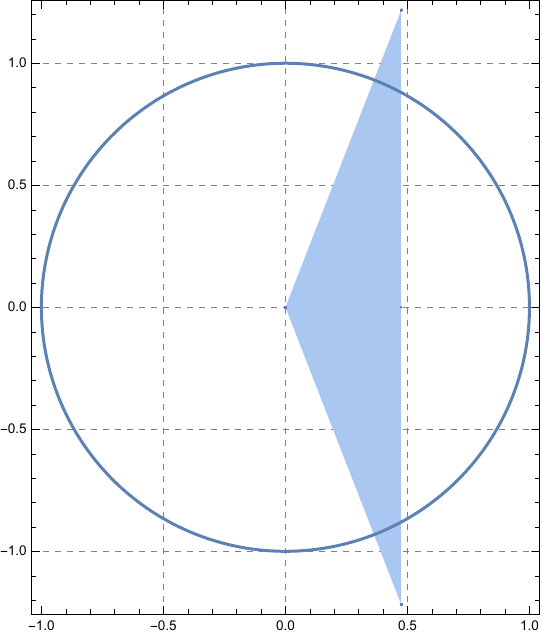}
            \phantom{xx} (b) Case $l=1/10$
        \end{minipage}

        \caption{$F(A^{-1}B)$ for $l=-1/10$ and $l=1/10$}

        \label{fig:4:erroresThetaMétodoEjemplo1}
    \end{figure}
\end{example}

\begin{example}
    \label{ejem:5:FisherKolmogorov2D}
    We consider now a variant of the Fisher-Kolmogorov equation in two space dimensions given by the delayed reaction-diffusion equation
    \begin{equation}
        \frac{\partial v}{\partial t} = \lambda \left( \frac{\partial^2 v}{\partial x^2} + \frac{\partial^2 v}{\partial y^2} \right) + \mu v(t-\tau,x,y)(1-v(t-\tau,x,y)), \quad t \geq 0, \quad 0 \leq x, y \leq 1,
        \label{eq:4:EDPR2}
    \end{equation}
    with the initial and boundary conditions
    \begin{equation*}
        \begin{cases}
        v(t,0,y)=v(t,1,y)=0,\\
        v(t,x,0)=v(t,x,1)=0,\\
        v(t,x) = \sin (\pi x) \sin(\pi y),\quad t \leq 0,
        \end{cases}
    \end{equation*}
    where $\lambda$ and $\mu$ are two positive parameters and $\tau$ is as before.

    On this occasion we define $\Delta x:= 1/M$, as well as the mesh points $x_0 = 0 < x_1 < \ldots < x_j = j \Delta x < \ldots < x_M=1$. For simplicity we will also denote $\Delta y = 1/M$ and $y_0 = 0 < y_1 < \ldots < y_j = j \Delta y < \ldots < y_M=1$. Then, denoting by $z^{i,j}$ the approximation of $v(t,x_i,y_j)$ and by $z^j$ the vector $(z^{1,j}, \ldots, z^{M-1,j})^T$, we obtain the MOL approach
    \begin{equation}
        z'(t) = A z(t) + g(z(t-\tau)),
        \label{eq:4:EDPR2_Adaptada}
    \end{equation}
    where now $z(t)=(z^1(t), \ldots, z^{M-1}(t))^T$,
    \begin{equation*}
        A =  \frac{\lambda}{\Delta y^2} \begin{pmatrix}
            -2 I + \Delta y ^2 L & I  & 0 & \ldots & 0\\
            I & -2 I + \Delta y ^2 L & I &\ldots & 0\\
            0 & I & -2 I + \Delta y ^2 L & \ldots & 0\\
            \vdots & \vdots & \vdots & \ddots & \vdots\\
            0 & 0 & \ldots  & I & -2 I + \Delta y ^2 L\\
        \end{pmatrix}, \quad
        L = \frac{1}{\Delta x^2}
        \begin{pmatrix}
            -2 & 1 & 0 & \ldots & 0\\
            1 & -2 & 1 & \ldots & 0\\
            0 & 1 & -2 & \ldots & 0\\
            \vdots & \vdots & \vdots & \ddots & \vdots\\
            0 & 0 & \ldots & 1 & -2\\
        \end{pmatrix},
    \end{equation*}
    $I$ is the identity matrix of dimension $M-1$
    and the \textit{j}th component of $g(z(t-\tau))$ is the vector
    \begin{equation*}
        g_j(z(t-\tau))^T = \left [ \mu z^{i,j}(t-\tau)(1-z^{i,j}(t-\tau)) \right ]_{i=1,\ldots,M-1}.
    \end{equation*}

    That is, we can also calculate a numerical solution of (\ref{eq:4:EDPR2}) by applying some $\theta$-method to the equation (\ref{eq:4:EDPR2_Adaptada}). For this reason, we will consider the $\theta$-method with $u=0$ and $\theta=1$ and will study for which parameters we can guarantee its unconditional stability when it is applied to (\ref{eq:4:EDPR2_Adaptada}).

    Firstly, for every $j=1, \ldots, M-1$ it is clear that each one of the components of the function $g_j$, say $g_{j,1}, \ldots, g_{j,M-1}$, can be understood as a differentiable function of a single variable that vanishes in $0$. Therefore, by the mean value theorem and through an abuse of notation, we can write $g_{j,i}(z) = g_{j,i}'(c_{i,j})$ for some $c_{i,j} \in \left(0,z^{i,j} \right)$.

    In other words, to study the stability of the $\theta$-method with $u=0$ and $\theta = 1$ for (\ref{eq:4:EDPR2_Adaptada}), we only have to study the stability of
    \begin{equation}
           z_{n+1} = z_n + h A \cdot z_{n+1} + h B_{n-m+1} \cdot z_{n-m+1},
           \label{eq:4:ThetaMetodoConB_Variable}
    \end{equation}
    where the diagonal matrix
    \begin{equation*}
        B_{n-m+1} := \begin{pmatrix}
            b_1 & 0 & \ldots & 0\\
            0 & b_2 & \ldots & 0\\
            \vdots & \vdots  & \ddots & \vdots\\
            0 & 0 & \ldots  & b_{M-1}\\
        \end{pmatrix}, \quad
        b_j := \begin{pmatrix}
            g_{j,1}'(c_{1,j}) & 0 & \ldots & 0\\
            0 & g_{j,2}'(c_{2,j}) & \ldots & 0\\
            \vdots & \vdots & \ddots & \vdots\\
            0 & 0 & \ldots & g'_{j,M-1}(c_{M-1,j})\\
        \end{pmatrix},
    \end{equation*}
    depends on the vector $z_{n-m+1}$, since each $c_{i,j}$ belongs to the interval $ \big(0,z_{n-m+1}^{i,j} \big)$. However, from Theorem \ref{thm:3:EstabilidadThetaMetodo} it can be deduced that if the inclusion $F\left( B_{n-m+1} A^{-1} \right) \subseteq D(0,1)$ is verified for all $n$, then the method (\ref{eq:4:ThetaMetodoConB_Variable}) is unconditionally stable. Therefore we only have to see for which parameters the above hypothesis is true.

    We assume that for all $n$ it is verified
    \begin{equation}
        \big |z_{n-m+1}^{i,j} \big | \leq 1, \quad i, j=1,\ldots,M-1.
        \label{eq:4:hipótesisCálculoParámetros}
    \end{equation}
    Then $-\mu < g_{j,i}'(c_{i,j}) < 3\mu$ and, since $B_{n-m+1}$ is always a hermitian matrix, we have that
    \begin{equation*}
        F(B_{n-m+1}) = \left[\min \limits_{i,j} (g_{j,i}'(c_{i,j})), \max \limits_{i,j} (g_{j,i}'(c_{i,j})) \right] \subseteq (-\mu, 3\mu).
    \end{equation*}

    On the other hand, the eigenvalues of $A$ are known, because $L$ is a Toeplitz matrix, which can be orthogonally diagonalisable as
     \begin{equation*}
     L = Q \diag(  \omega_1, \ldots, \omega_{M-1} ) Q^T = Q \Lambda Q^T,
        \end{equation*}
    where $\omega_k := -4 \lambda M^2 \sin^2 \left ( \frac{k \pi}{2 M} \right )$.

    Thus, $ A = (I \otimes Q) \tilde{A} (I \otimes Q^T)$, where $\tilde{A}$ is the block-Toeplitz matrix obtained by replacing $L$ in $A$ by $\Lambda$.  Hence,
     \begin{equation*}
     A = (I  \otimes  Q) K ( L_1  \oplus  \ldots \oplus L_{M-1} ) K (I \otimes  Q^T),
        \end{equation*}
    where $K$ is the commutation matrix, and
         \begin{equation*}
       L_i =
        \begin{pmatrix}
           \omega_i  & \frac{1}{\Delta x^2} & 0 & \ldots & 0\\
             \frac{1}{\Delta x^2} & \omega_i &  \frac{1}{\Delta x^2} & \ldots & 0\\
            0 &  \frac{1}{\Delta x^2} & \omega_i & \ldots & 0\\
            \vdots & \vdots & \vdots & \ddots & \vdots\\
            0 & 0 & \ldots &  \frac{1}{\Delta x^2} & \omega_i\\
        \end{pmatrix},
    \end{equation*}
    consequently, the eigenvalues of $A$ are of the form $\omega_i + \omega_j$, $i,j=1,\ldots,M-1$.
    Therefore
    \begin{equation*}
        \sigma(A) \subseteq [2\omega_{M-1}, 2\omega_1] = \left [ -8 \lambda M^2 \sin^2 \left ( \frac{(M-1) \pi}{2 M} \right ), -8 \lambda M^2 \sin^2 \left ( \frac{\pi}{2 M} \right ) \right].
    \end{equation*}
    But what is more, the matrix $A$ is negative definite, so it also is hermitian. That is, it is verified that the field of values
    $F(A) =  Co (\sigma(A))$ and, as a result, that $F(A)$ is also included in the above interval.

    As a result, if (\ref{eq:4:hipótesisCálculoParámetros}) is verified, then (in a similar way as in \cite[Prop.~6.1]{SeiboldPractice}):
    \begin{equation*}
    F(  A^{-1/2} B_{n-m+1} A^{-1/2} ) = \{ <v, B_{n-m+1}v> : <v,  Av> = 1, v \in \mathbb{V} \} =
    \{ <v, B_{n-m+1} v> : <v,  v> = 1/ l, l \in F(A),   v \in \mathbb{V} \}  \end{equation*}
     \begin{equation*}
     \subseteq \left[-\frac{3\mu}{8 \lambda M^2 \sin^2 \left ( \frac{\pi}{2 M} \right )}, \frac{\mu}{8 \lambda M^2 \sin^2 \left ( \frac{\pi}{2 M} \right )} \right].
    \end{equation*}

    But furthermore, for every $t \leq 0$ and $i=1,\ldots,M-1$, we have that $\big |z^i(t) \big | \leq 1$. Hence, it can also be proved (by induction) that if
    \begin{equation*}
        \Bigg |-\frac{3\mu}{8 \lambda M^2 \sin^2 \left ( \frac{\pi}{2 M} \right )} \Bigg| < 1,
    \end{equation*}
    the hypothesis (\ref{eq:4:hipótesisCálculoParámetros}) is always true, and, consequently, whenever the parameters of (\ref{eq:4:EDPR2_Adaptada}) verify the condition
    \begin{equation}
        \lambda > \frac{3\mu}{8 M^2 \sin^2 \left ( \frac{\pi}{2 M} \right )},
        \label{eq:4:condEstabilidadIncondicional}
    \end{equation}
    the method (\ref{eq:4:ThetaMetodoConB_Variable}) is unconditionally stable.

\begin{figure}[htb]
    \centering

    \begin{minipage}[t]{0.45\textwidth}
        \centering
        \includegraphics[width=\textwidth,height=5cm]{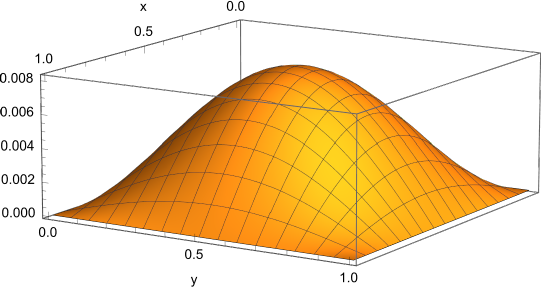}
         (a) For $t=3$
    \end{minipage}
    \hfill
    \begin{minipage}[t]{0.52\textwidth}
        \centering
        \includegraphics[width=\textwidth,height=5cm]{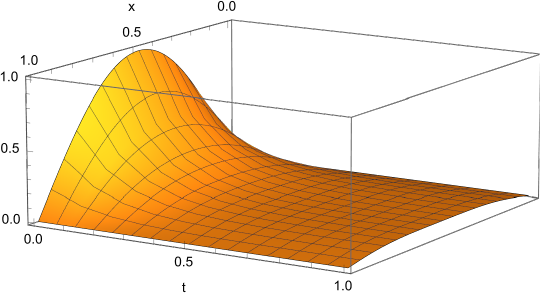}
        (b) For $y=0.5$
    \end{minipage}

    \caption{Numerical solution of (\ref{eq:4:EDPR2}) given by the $\theta$-method with $u=0$ and $\theta=1$}

    \label{fig:4:soluciónThetaMétodoEjemplo2}
\end{figure}

    This way, we know that if we set, for example, $\lambda=1/2$, $\mu=3$ and $M=100$, the $\theta$-method with $u=0$ and $\theta=1$ will always be stable, regardless of the step size considered. Fig. \ref{fig:4:soluciónThetaMétodoEjemplo2} shows the numerical solution given by this $\theta$-method if we use $\tau=1$ and $m=10$, besides the above parameters. Note that the problem is symmetric in the sense that we obtain the same solution if we interchange the spatial variables in (\ref{eq:4:EDPR2}).
\end{example}

\section{Conclusions and future challenges}

The main focus of this paper was to show that we can also apply to problems with delay the
relatively new theory used in \cite{SeiboldTheory,SeiboldPractice} to study the stability of certain numerical methods for ODEs, since, in this way, we have been able to present a possible approach to study the asymptotic behavior of DDEs and PDDEs more general than the one we have always worked with.

In this sense, several remarks need to be noted. With this theory we have not only presented a sufficient condition for the asymptotic stability of any $\theta$-method for the DDE (\ref{eq:1:EDR_EstabilidadThetaMétodo}) in Proposition \ref{prop:3:EstabilidadThetaMetodo}, but we have also shown another condition that guarantees the unconditional stability of any $\theta$-method with $u=0$ and $\theta>1/2$ for this same equation with Theorem \ref{thm:3:EstabilidadThetaMetodo}, which has allowed us to work more easily with the numerical examples for PDDEs proposed here. In addition, we have also characterized the regions $D_y$ given by $u=0$ and $\theta=1$ (and the relation between them) to facilitate the study in those cases where we cannot guarantee unconditional stability. Finally, we have also seen how to apply these results to certain parabolic problems given by PDDEs with one diffusion term and one delayed term.

However, there are still several issues that remain which require further analysis. First of all, in this paper we have only focused on $\theta$-methods for DDEs, so it would be ideal to extend the theory presented here to a wider range of methods for problems with delay. Moreover, for simplicity, when presenting the numerical examples for PDDEs we have only considered the $\theta$-method with $u=0$ and $\theta=1$, therefore another future challenge could be to find equivalent results for other values of these two parameters.
The latter would also be appropriate when studying unconditional stability or
characterizing the regions $D_y$.
On the other hand, all partial delay differential equations considered follow the same structure (one diffusion term and one delayed term), so it would also be interesting to try to apply the theory discussed here to other types of problems with delay (including more systems of PDDEs or problems in more dimensions).

\section*{Acknowlgedgments}
   The authors would like to thank  Karel in't Hout and  David Shirokoff for useful discussions and comments.   This research was funded by the Spanish Ministerio de Ciencia e Innovación (MCIN) with funding from the European Union NextGenerationEU (PRTRC17.I1) and the Consejería de Educación, Junta de Castilla y Le\'on, through QCAYLE project, and also by Fundaci\'on Sol\'orzano through the project FS/5-2022.
\vspace{0.1cm}

\noindent
\small
\textbf{Data Availability} Data sharing is not applicable to this article as no datasets were generated or analyzed during the current study.

\noindent
\textbf{Code Availability} The codes employed for the current study are available upon request.

\section*{Declarations}
\small
\textbf{Conflict of interest}
The authors declare that they have no conflict of interest.


\begin{thebibliography}{}

\bibitem{BoRi00}
G.~A. Bocharov, F.~A. Rihan, Numerical modelling in biosciences using delay
  differential equations, Journal of Computational and Applied Mathematics
  125~(1) (2000) 183--199, numerical Analysis 2000. Vol. VI: Ordinary
  Differential Equations and Integral Equations.
\newblock \href
  {http://dx.doi.org/https://doi.org/10.1016/S0377-0427(00)00468-4}
  {\path{doi:https://doi.org/10.1016/S0377-0427(00)00468-4}}.

\bibitem{macdonald2008biological}
N.~MacDonald, N.~MacDonald, C.~Cannings, F.~Hoppensteadt,
  \href{https://books.google.es/books?id=eRksych27yQC}{Biological Delay
  Systems: Linear Stability Theory}, Cambridge Studies in Mathematical Biology,
  Cambridge University Press, 2008.

\bibitem{MaGl77}
M.~Mackey, L.~Glass, Oscillation and chaos in physiological control systems,
  Science (New York, N.Y.) 197 (1977) 287--9.
\newblock \href {https://doi.org/10.1126/science.267326}
  {\path{doi:10.1126/science.267326}}.

\bibitem{TaMaBe00}
Y.~Takeuchi, W.~Ma, E.~Beretta, Global asymptotic properties of a delay sir
  epidemic model with finite incubation times, Nonlinear Anal. 42~(6) (2000)
  931–947.
\newblock \href {http://dx.doi.org/10.1016/S0362-546X(99)00138-8}
  {\path{doi:10.1016/S0362-546X(99)00138-8}}.

\bibitem{Jackiewicz1987/88}
B.~A. Z.~M. Jackiewicz, Z., \href{http://eudml.org/doc/133256}{Stability
  analysis of one-step methods for neutral delay-differential equations.},
  Numerische Mathematik 52~(6) (1987/88) 605--620.


Numerical Modeling by Delay and Volterra Functional Differential Equations

\bibitem{Baker06}
C.~Baker, G.~Bocharov, A.~Filiz, N.~Ford, C.~Paul, F.~Rihan, A.~Tang,
  R.~Thomas, H.~Tian, D.~Wille, Numerical Modeling by Delay and Volterra Functional Differential Equations, 2006.

\bibitem{KotoRK}
Koto, T.: Stability of IMEX Runge-Kutta methods for delay differential equations. Journal of Computational and Applied Mathematics \textbf{211}, 201-212 (2008). \url{https://doi.org/10.1016/j.cam.2006.11.011}

\bibitem{KotoMultistep}
Koto, T.: Stability of implicit-explicit linear multistep methods for ordinary and delay differential equations. Frontiers of Mathematics in China \textbf{4}, 113-129 (2009). \url{https://doi.org/10.1007/s11464-009-0005-9}


\bibitem{baker_buckwar_2000}
C.~T.~H. Baker, E.~Buckwar, Numerical analysis of explicit one-step methods for
  stochastic delay differential equations, LMS Journal of Computation and
  Mathematics 3 (2000) 315–335.
\newblock \href {http://dx.doi.org/10.1112/S1461157000000322}
  {\path{doi:10.1112/S1461157000000322}}.

\bibitem{KUCHLER2000189}
U.~Küchler, E.~Platen, Strong discrete time approximation of stochastic
  differential equations with time delay, Mathematics and Computers in
  Simulation 54~(1) (2000) 189--205.
\newblock \href
  {http://dx.doi.org/https://doi.org/10.1016/S0378-4754(00)00224-X}
  {\path{doi:https://doi.org/10.1016/S0378-4754(00)00224-X}}.



\bibitem{Calvo}
Calvo, M., Grande, T.: On the asymptotic stability of $\theta$-methods for delay differential equations. Numerische Mathematik \textbf{54}, 257-269 (1988). \url{https://doi.org/10.1007/BF01396761}


\bibitem{Rihan}
Rihan, F.A.: Delay differential equations and applications to biology. Springer, Heidelberg (2021). \url{https://doi.org/10.1007/978-981-16-0626-7}

\bibitem{DAUTILIA20202067}
M.~C. D’Autilia, I.~Sgura, V.~Simoncini, Matrix-oriented discretization
  methods for reaction–diffusion {PDE}s: Comparisons and applications,
  Computers and Mathematics with Applications 79~(7) (2020) 2067--2085, advanced
  Computational methods for PDEs.
\newblock \href {http://dx.doi.org/https://doi.org/10.1016/j.camwa.2019.10.020}
  {\path{doi:https://doi.org/10.1016/j.camwa.2019.10.020}}.

\bibitem{Cox2002ExponentialTD}
S.~M. Cox, P.~C. Matthews, Exponential time differencing for stiff systems,
  Journal of Computational Physics 176 (2002) 430--455.

\bibitem{VIGOAGUIAR200780}
J.~Vigo-Aguiar, J.~Martín-Vaquero, Exponential fitting {BDF} algorithms and
  their properties, Applied Mathematics and Computation 190~(1) (2007) 80--110.
\newblock \href {http://dx.doi.org/https://doi.org/10.1016/j.amc.2007.01.008}
  {\path{doi:https://doi.org/10.1016/j.amc.2007.01.008}}.

\bibitem{VAMVWa07}
J.~Vigo-Aguiar, J.~Martín-Vaquero, B.~A. Wade, Adapted {BDF} algorithms applied
  to parabolic problems, Numerical Methods for Partial Differential Equations
  23~(2) (2007) 350--365.
\newblock \href {http://dx.doi.org/https://doi.org/10.1002/num.20180}
  {\path{doi:https://doi.org/10.1002/num.20180}}.

\bibitem{Higueras_2015}
I.~Higueras, T.~Rold{\'{a}}n,
 Construction of additive
  semi-implicit {R}unge--{K}utta methods with low-storage requirements,
  Journal of Scientific Computing 67~(3) (2015) 1019--1042.
\newblock \href {https://doi.org/10.1007/s10915-015-0116-2}
  {\path{doi:10.1007/s10915-015-0116-2}}.

\bibitem{SeiboldTheory}
Rosales, R.R., Seibold, B., Shirokoff, D., Zhou, D.: Unconditional Stability for Multistep ImEx Schemes: Theory. SIAM Journal of Numerical Analysis \textbf{55}, 2336-2360 (2017). \url{https://doi.org/10.1137/16M1094324}

\bibitem{SeiboldPractice}
Seibold, B., Shirokoff, D., Zhou, D.: Unconditional stability for multistep ImEx schemes: Practice. Journal of Computational Physics \textbf{376}, 295-321 (2019). \url{https://doi.org/10.1016/j.jcp.2018.09.044}

\bibitem{kj1994stability}
K.~in't Hout, The stability of $\theta$-methods for systems of delay
  differential equations, Ann. Numer. Math. 1 (1994) 323–--334.

\bibitem{HornJohnson}
Horn, R.A., Johnson, C.R.: Topics in Matrix Analysis. Cambridge University Press, Cambridge (1991)

\bibitem{Johnson}
Johnson, C.R.: Numerical determination of the field of values of a general complex matrix. SIAM Journal of Numerical Analysis \textbf{15}, 595-602 (1978). \url{https://doi.org/10.1137/0715039}

\bibitem{Bellen}
Bellen, A., Zennaro, M.: Numerical Methods for Delay Differential Equations. Oxford University Press, Oxford (2013)

\end{thebibliography}
\end{document}